\documentclass[11pt,a4paper]{amsart}
\usepackage{graphicx}
\usepackage[labelsep=none]{caption}
\usepackage{amsfonts}
\usepackage{amsthm}
\newtheorem{theorem}{Theorem}[section]
\newtheorem{lemma}[theorem]{Lemma}

\title{On some impossible disentanglement puzzles}
\author{Fernando Galve Mauricio}
\email{fgalvem@yahoo.es}
\dedicatory{Dedicado a la memoria de mi padre \\
To my father, in loving memory}
\begin{document}
\begin{abstract}
In this article we prove the impossibility of some disentanglement puzzles, first building mathematical models that reflect the essential characteristics of these puzzles.
\end{abstract}
\maketitle

\section*{Introduction}
A disentanglement puzzle is a mechanical puzzle consisting of two or more pieces made of metal, wood or string that one is required to separate subject to the constraints that the solver cannot bend the wires or break the loops of rope. The pieces aren't linked together, so if they were arbitrarily deformable, solving the puzzle would usually be trivial \cite{MatthewHorak2006}.

In the literature there are some proofs of impossibility of puzzles such us Steward Coffin's one, that amounts to remove a loop of string wrapped around a rigid wire \cite{IntaBertuccioni2003}, or plenty of similar impossible puzzles studied in \cite{PaulMelvin2004}. But the point there is that these are not disentanglement puzzles in the sense we employ in this article, because of the (nontrivial) fact that the loop and the wire are topologically linked.

In this paper we are interested in the impossibilities arising not from topology, but from the concrete sizes of the rigid pieces and/or the lengths of the loops of rope.

\section{The simplest disentanglement puzzle}
\begin{figure}
\includegraphics[scale=0.15]{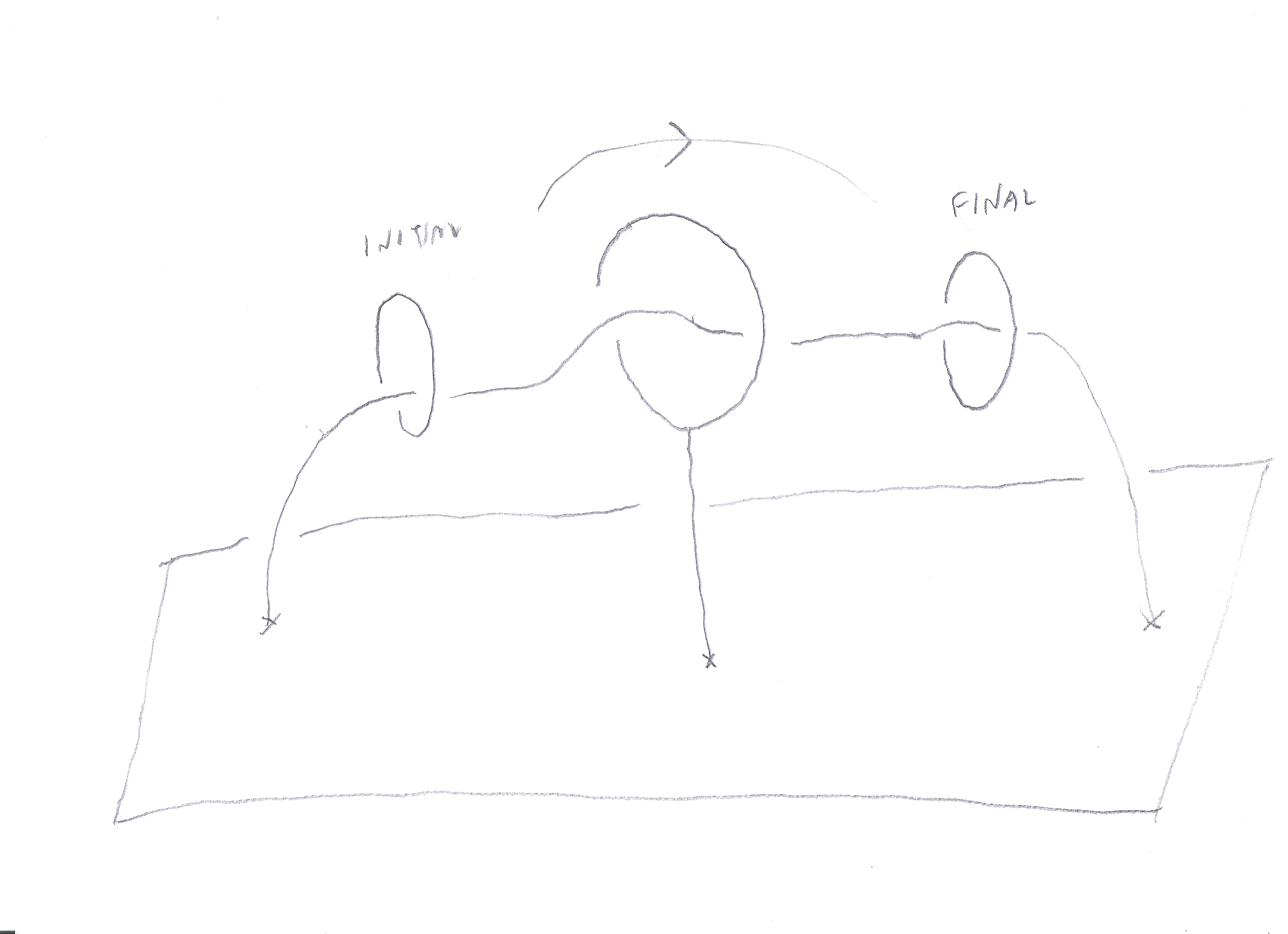}
\caption{}
\label{ModelAPuzzle}
\end{figure}

Here we will study the puzzle consisting of a pole of wire ended in a circle, fixed on a plane, throug which passes a rope that is fixed on both ends to the plane, and the third piece is another wire or wooden hoop that surrounds the rope. The goal of the puzzle is to pass the hoop through the wire circle such that it ends at the other side of it (see figure~\ref{ModelAPuzzle}) and everything is realized in one half-space of the plane.
We will model this puzzle mathematically and prove that it is solvable if and only if the radius of the hoop is less than the radius of the circle upon the pole.

\subsection{The \textsl{Model A} puzzle}
Here is the mathematical model of the preceding puzzle which will be called \textsl{Model A} puzzle.

The table to which are attached the pole and the string is the plane $z=0$ of $\mathbb{R}^3$ with the usual x, y and z coordinates. The pole and its ending circle are the segment from $(0,0,0)$ to $(0,0,1)$ and the circle $\{x^2+y^2+(z-2)^2=1\} \cap \{y=0\}$. The string of rope is the poligonal line passing through $(0,-1,0)$, $(0,-1,2)$, $(0,1,2)$ and $(0,1,0)$. Finally, the moving hoop is initially $\{x^2+(y+\frac{1}{2})^2+(z-2)^2=r^2\} \cap \{y=-\frac{1}{2}\}$ for some fixed $1 \leq r < 2$.

The goal is to pass from this configuration to the exact same, but with the moving hoop in the position $\{x^2+(y-\frac{1}{2})^2+(z-2)^2=r^2\} \cap \{y=\frac{1}{2}\}$ allowing moving the hoop and performing homotopies of the rope that do not touch each other neither touch the plane $z=0$ nor the pole with its ending circumference.

\subsection{Two easy lemmas}
To prove that solving \textsl{Model A} puzzle is impossible we will invoke two easy lemmas whose proof are left to the reader.

\begin{lemma}\label{lemmaCircle}
If in $\mathbb{R}^2$ we have a fixed circle of radius $r \geq 1$ and a segment is moving continuously such as its length variable is less or equals than 1, and such that intersects the circle in two points but neither of its endpoints is allowed to cross the circumference, then we can continously select one of the sides of the segment such that it contains less than half the circle.
\end{lemma}
\proof Left to the reader.\qed

\begin{lemma}\label{lemmaLollipop}
If in $\mathbb{R}^2$ we have the figure formed by the pole $(0,0)$ to $(0,1)$ and the circle $x^2+(y-2)^2=1$ being cut from a continuously moving segment such that its ends do not touch the circle nor the pole, nor any end of the segment is inside the circle, nor the segment passes through $(0,0)$ then we can continuously assign one of the sides of the segment such that the portion of the lollipop cutted by the segment on this side does not contain the point $(0,0)$.
\end{lemma}
\proof Left to the reader.\qed

\subsection{Impossibility of \textsl{Model A} puzzle}
\begin{theorem}
It is impossible to solve the \textsl{Model A} puzzle.
\end{theorem}
\proof It suffices to prove that we can construct continously varying spanning disks of the hoop and of the circle upon the pole in such a way that they are simple flat circle when they are apart, because then the string of rope traverses those disks in order giving origin to a $\mathbb{Z}*\mathbb{Z}$ invariant.

The spanning disks will be the flat circles when they do not touch each other nor the pole. When the flat disk spanning the hoop cuts the flat disk spanning the circle upon the pole, we use lemma~\ref{lemmaCircle} to replace a neighborhood of the intersection segment with a continously varying hood, so that both surfaces avoid each other. And when the flat disk spanning the hoop is cutted by the flat disk spanning the circle upon the pole and/or the pole from the inside, we use lemma~\ref{lemmaLollipop} to build a continuously varying hood such that the deformed circle avoids the lollipop. \qed

\section{The \textsl{Model B} puzzle}
\begin{figure}
\includegraphics[scale=0.15]{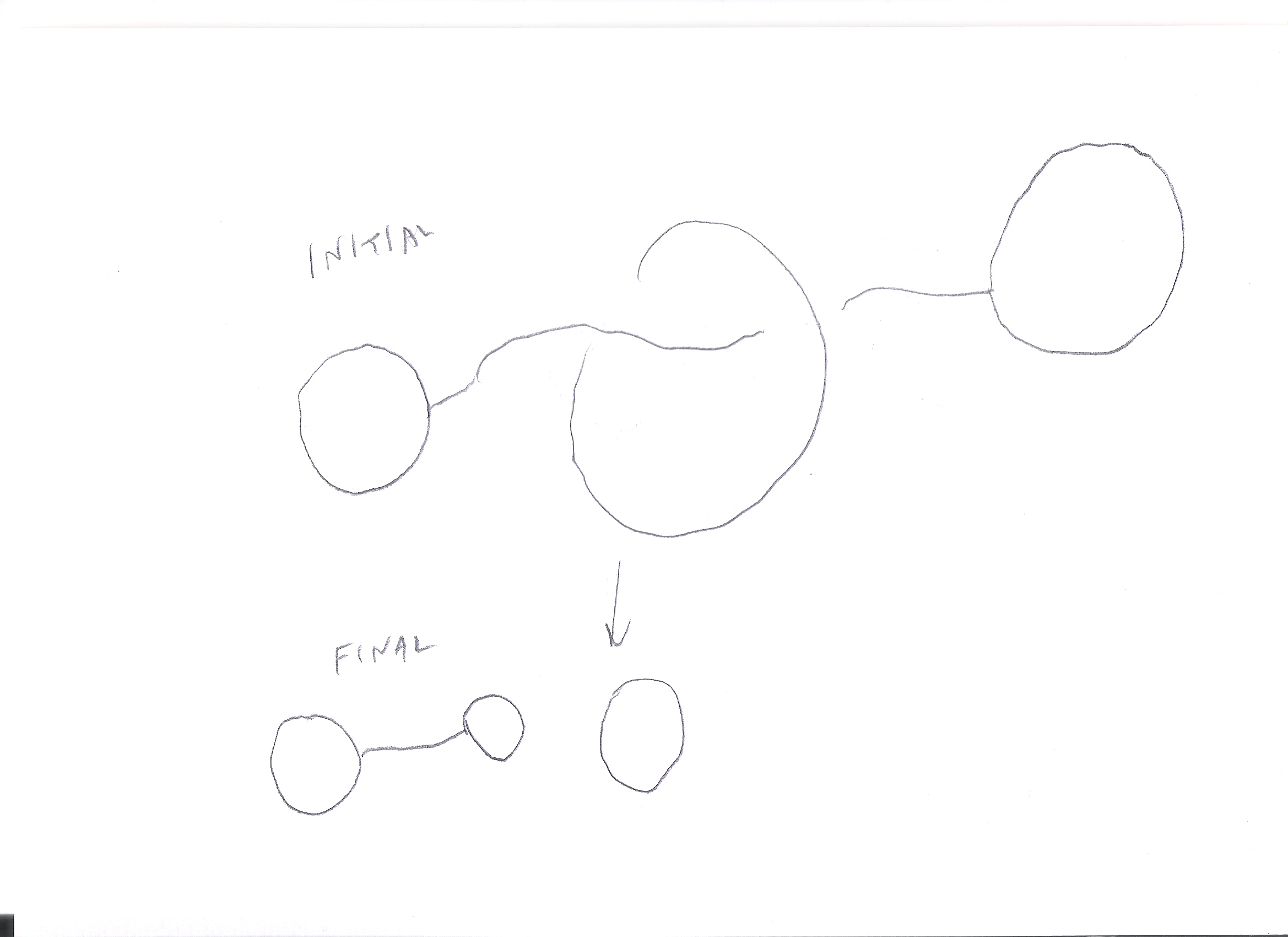}
\caption{}
\label{ModelBPuzzle}
\end{figure}

The \textsl{Model B} puzzle is an analogous puzzle consisting of a string of rope ended with a circumference at each side, and a rigid hoop in between, and the goal of the puzzle is to remove the middle hoop. We'll prove that if the radius of both circumferences are greater or equal to the radius of the rigid hoop, then the puzzle becomes impossible (see figure~\ref{ModelBPuzzle}).

\begin{theorem}
\textsl{Model B} puzzle is impossible even if we allow arbitrary homotopies of the rope, and the rope to cut the circles at its ends, and even the circles to cut each other, as long as we forbid the rope and the circles at its end to touch the middle hoop.
\end{theorem}
\proof Without loss of generality we may assume that the middle hoop is fixed and only the circles and the rope joining them are moving.

The universal covering of the middle hoop complement is divided into infinitely many equal parts by the lifts of its spaning circle, in such a way that when we lift the circles and the rope joining then to this universal covering, each circle lies at a different part.

Now, it suffices to prove that when we lift the homotopy to this universal covering, at the end each circle ends at a different part; and for this it suffices to prove that at each moment of the homotopy, each lifted circle has at least one point in its initial part of the covering, but this is again a trivial consequence of lemma~\ref{lemmaCircle} above. \qed

\section{The \textsl{Model C} puzzle}
\begin{figure}
\includegraphics[scale=0.10]{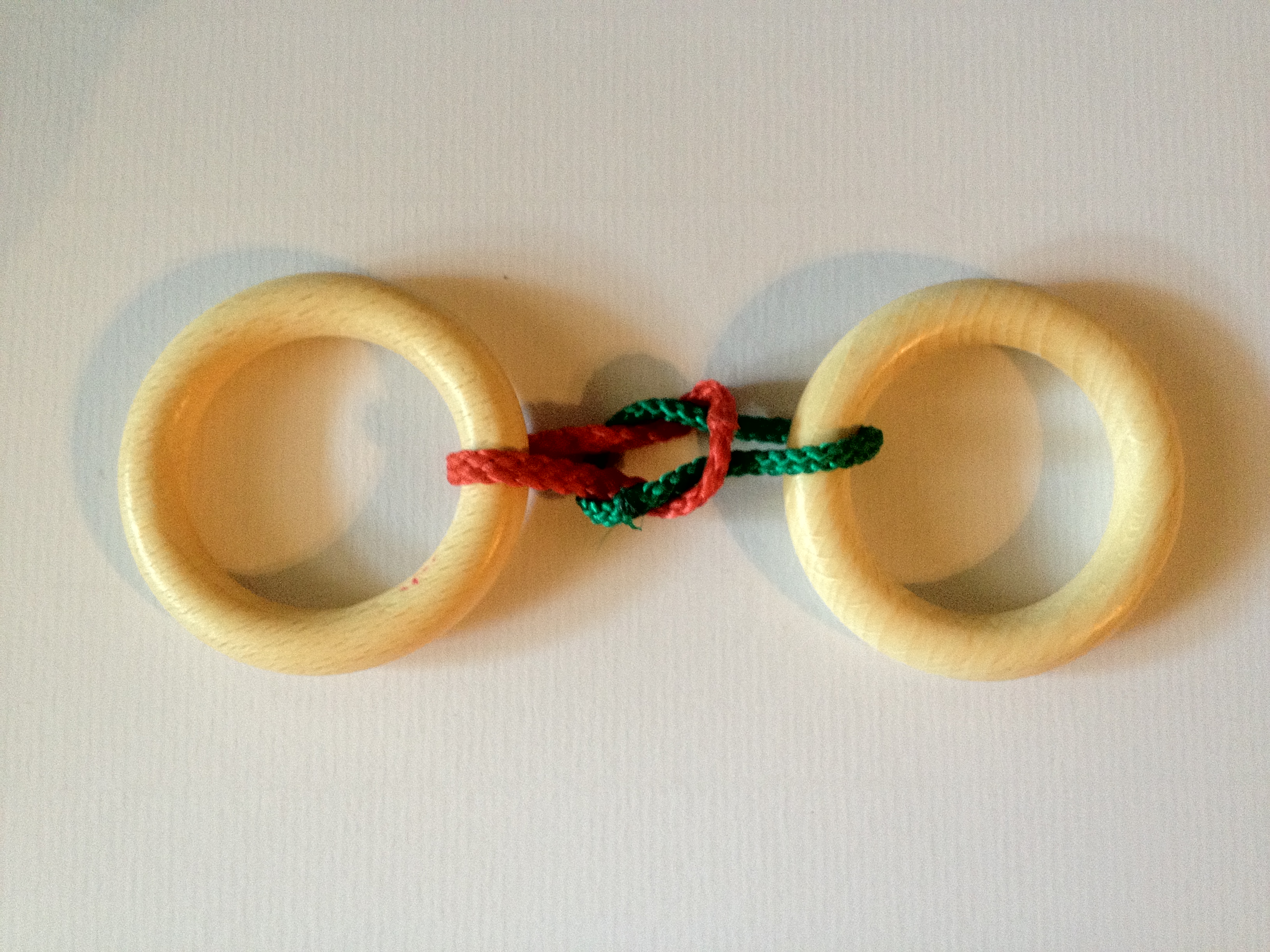}
\caption{}
\label{ModelCPuzzle}
\end{figure}

The \textsl{Model C} puzzle is formed by two Hopf links formed each by a wooden (solid) hoop and a loop of rope linked to it. Both hoops are the same radius and both ropes are the same length, though usually shorter than the hoops. The links are tangled as shown in figure~\ref{ModelCPuzzle} and the goal is to separate them.
In this case, if we allow arbitrary homotopies of the ropes, the puzzle is easily solved; so instead of homotopies of the ropes we only allow $\Delta$-moves. I mean, each rope is modeled by a simple closed polygon, and each $\Delta$-move consists of taking a solid triangle than only has one side in common with the polygon and does not touch the other polygon neither any of the hoops, and replacing this side by the two other sides (plus the inverse operation).

\subsection{An auxiliary lemma}
\ 

In the following lemma we work in $S^3=\mathbb{R}^3 \cup \{\infty\}$.
\begin{lemma}\label{lemmaThreeArcs}
If we make an ambient isotopy of two unlinked circles (retaining its circular shape) plus a simultaneous isotopy of a sphere (also keeping it always spherical) such that at every moment the sphere cuts an arc of every circle, then there is an extending isotopy of a figure formed by the two circles joined by an arc of circumference always in the outer side of the moving sphere. Furthermore, the outer side of the sphere minus the three arcs of circle (from the two initial cirles plus the joining one) retracts by isotopic deformation to a neighborhood of the joining arc of circle.
\end{lemma}
\proof Select one of the endpoints of the arcs of circle in which the sphere cuts one of the two unlinked circles, continously moving with the isotopies, and make an spherical inversion with center that point. The unbounded component of the sphere and the two arcs of circle transform into a half space (limited by a plane), plus a straight line ray, plus an arc of circle with its ends in the limiting plane. Join the straight line to the arc of circle by the horizontal segment from the point of the arc that is farthest from the plane, to the ray of line. Inverting again, this segment transform to the requested arc of circle joining the two isotoping circles.

The rest of the proof (the existence of the contraction to a neighborhood of the arc joining both circles) is easy and is left to the reader.
\qed

\subsection{The impossibility proof}
\begin{figure}
\includegraphics[scale=0.15]{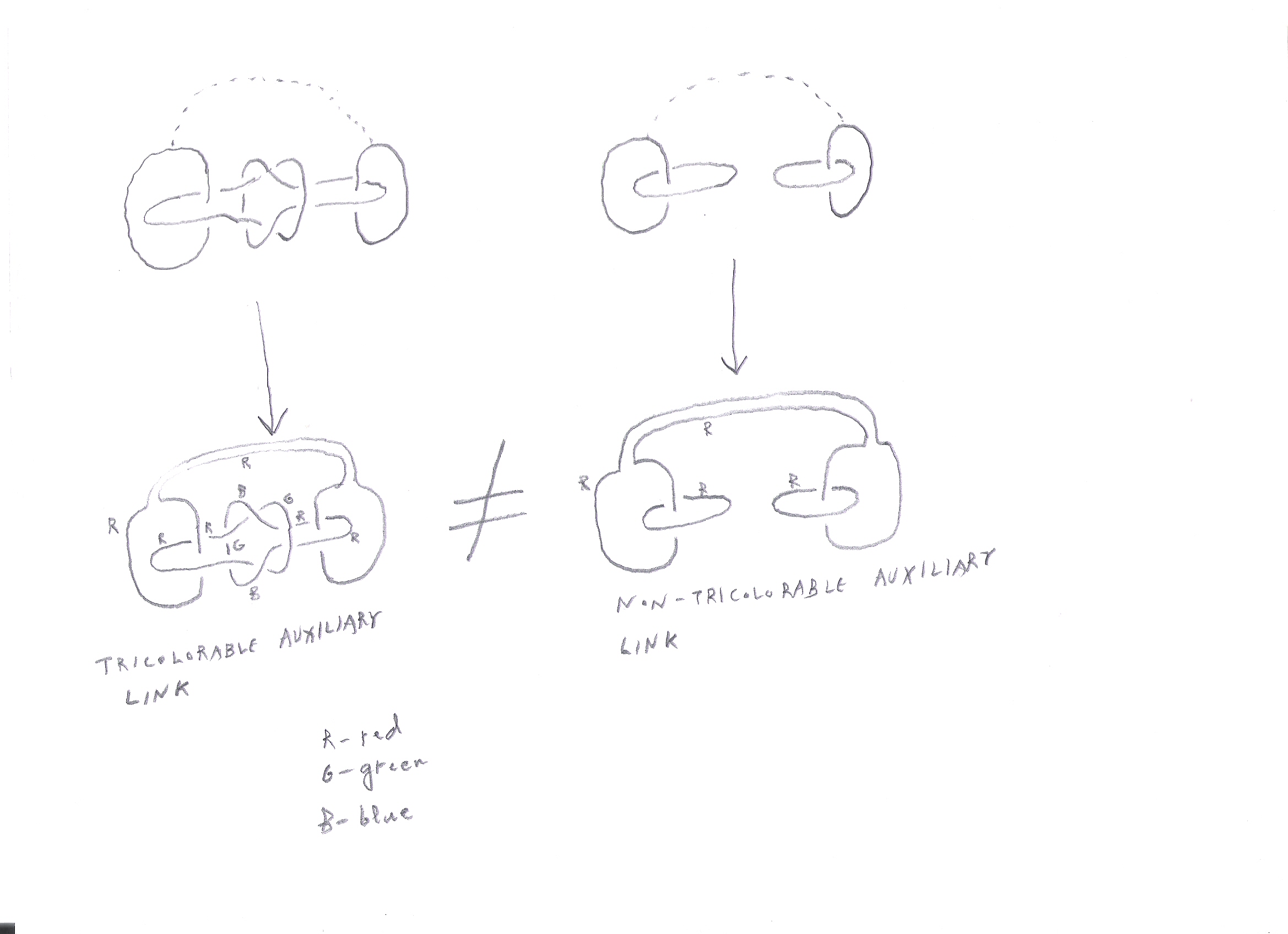}
\caption{}
\label{ModelCPuzzleAuxiliaryLinks}
\end{figure}

\begin{theorem}
\textsl{Model C} puzzle is impossible if we maintain the lengths of the ropes less than or equal to the radius of the hoops.
\end{theorem}
\proof We proceed by contradiction, assuming that it is possible to separate the pieces of the puzzle. Using a band sum of the Hopf links as seen in figure~\ref{ModelCPuzzleAuxiliaryLinks} we construct an auxiliary tricolorable link; and starting from two untangled Hopf links we construct a second auxiliary non-tricolorable link, also in figure~\ref{ModelCPuzzleAuxiliaryLinks}. We shall construct an ambient isotopy of both links, hence the contradiction.

We can construct a moving auxiliary sphere of the same radius of the hoops with center one continously varying point in one of the ropes. Initially the sphere contains both ropes in its bounded component, and cuts an arc of every hoop, so we can join the two hoops by an arc of circle in the unbounded component of the sphere, such as we did in lemma~\ref{lemmaThreeArcs}.

Now, as we isotope the pieces of the puzzle to separate them, there is a first moment in which one of the ropes touches the sphere. In this moment we take the concentric sphere of radius half the radius of the hoops. This half sized sphere contains one of the ropes in its interior and lefts the other one in its exterior.

Now consider two cases.

If the last sphere cuts only one of the hoops, then its linked rope is in the sphere interior and can be isotoped to a curve arbitrarily close to the arc of hoop inside the sphere.

If the half sized sphere cuts both hoops, then is its exterior rope that can be isotoped to a curve arbitrarily close to the arc of circle joining the two hoops, by the deformation statement of lemma~\ref{lemmaThreeArcs}.

In both cases we have now one of the ropes shrunken around an arc of circle. Now we just make the inverse isotopy, restoring the puzzle to its initial position, but maintaining the shrunken rope near its arc of circle, then isotope the other rope around its linked hoop, make again the isotopy forward, but maintaining the second rope around its hoop, undo the shrink of the rope, and finally do the backward isotopy and we end with both ropes apart each other, without ever touching the arc of circle connecting the two hoops.
\qed

Two final remarks with respect to \textsl{Model C}.
First, we conjecture the optimal length of the ropes, the limit between possible and impossible puzzles, to be $\frac{8}{3}$ the radius of the hoops.

The second remark is that using similar arguments, we can prove that the \textsl{Quattro} puzzle discussed in \cite{MatthewHorak2006} is impossible if the lengths of the ropes is maintained under one half the radius of the hoops.

\bibliography{OnSomeImpossibleDisentanglementPuzzles}
\bibliographystyle{ieeetr}

\end{document}